\pdfoutput=1
\documentclass{amsart}
\usepackage{graphicx}
\usepackage{hyperref}

\theoremstyle{plain}
\newtheorem{theorem}{Theorem}

\newtheorem{lemma}{Lemma}
\newtheorem{fact}{Fact}

\theoremstyle{definition}

\theoremstyle{remark}

\newcommand{\fonep}{
\,\begin{picture}(28,8)
\put(5,3){\circle{10}}
\put(25,3){\circle{10}}
\put(8.5,0){\vector(1,0){12}}
\put(20.5,6){\vector(-1,0){12}}
\put(10,2.5){\vector(-1,0){10}}
\put(2,3.8){\tiny $+$}
\end{picture}~
}

\newcommand{\fonem}{
\,\begin{picture}(28,8)
\put(5,3){\circle{10}}
\put(25,3){\circle{10}}
\put(8.5,0){\vector(1,0){12}}
\put(20.5,6){\vector(-1,0){12}}
\put(10,2.5){\vector(-1,0){10}}
\put(2,3.8){\tiny $-$}
\end{picture}~
}

\begin{document}
\title[Type II Reidemeister moves of links]{On Type II Reidemeister moves of links}
\author{Noboru Ito}
\address{National Institute of Technology, Ibaraki College,  312-8508, Japan}
\email{nito@gm.ibaraki-ct.ac.jp}
\date{January 17, 2022}
\maketitle

\begin{abstract}
\"{O}stlund (2001) showed that all planar isotopy invariants of generic plane curves that are unchanged under cusp moves and triple point moves, and of finite degree (in self-tangency moves) are trivial.   Here the term ``of finite degree" means Arnold-Vassiliev type. 
It implies the conjecture, which was often called \"{O}stlund conjecture:  
``Types I and III Reidemeister moves are sufficient to describe a homotopy from any generic immersion from the circle into the plain to the standard embedding of the circle".  Although counterexamples are known nowadays, there had been no (easy computable) function that detects the difference between the counterexample and the standard embedding on the plain.  However, we introduce a  desired function (Gauss diagram formula) is found for the two-component case.     
\end{abstract}


\section{Introduction}
Vassiliev \cite{Vassiliev1990} introduced knot invariants by using a stratification of the knot space by singularities.  Later Goussarov \cite{GoussarovPolyakViro2000} showed that every Vassiliev invariant is presented by a well-known notion, Gauss diagram  formula  \cite{GoussarovPolyakViro2000}.  

\"{O}stlund \cite{Ostlund2001PhD}  studied Vassiliev invariants using Gauss diagram formulas and observed that Gauss diagram formulas are always filtered by three types of singularities: cusps, self-tangencies, and triple points.  In fact, so are Arnold's  basic invariants of plane curves \cite{Arnold1994}.  Although Gauss diagram formulas  should have been already known as a useful tool for knots and plane curves (e.g., \cite{PolyakViro1994, Polyak1998}) then, \"{O}stlund, in the PhD thesis \cite{Ostlund2001PhD},  further unified two notions:  Gauss diagram formulas for knots and those of plane curves.   \"{O}stlund defined a notion of ``\emph{knot diagram invariant}" as follows: ``{\it We shall call a function of knot diagrams that is unchanged by planar isotopy, but not necessarily by Reidemeister moves, a knot diagram invariant}".    Let us  recall \"{O}stlund's results (Facts~\ref{fact:ostMove} and \ref{fact:ostInv}); the $n$th Reidemeister move is called $\Omega_n$-moves ($n=1, 2, 3$) here.  
\begin{fact}[{\cite[Chapter~IV, Theorem~2]{Ostlund2001PhD}}]\label{fact:ostMove}
Let $v$ be a knot diagram invariant that is unchanged under $\Omega_1$- and $\Omega_3$-moves, and of finite degree (in $\Omega_2$).  Then $v$ is a knot invariant.  
\end{fact}
\begin{fact}[{\cite[Chapter~IV, Corollary~1]{Ostlund2001PhD}}]\label{fact:ostInv}
All planar isotopy invariants of generic plane curves that are unchanged under cusp moves and triple point moves, and of finite degree (in self-tangency moves) are trivial.  
\end{fact}
That is, for \emph{knots}, any Gauss diagram formula (i.e., Vassiliev-type knot diagram invariant) does not detect the independence of $\Omega_2$-move or self-tangency move.  However, for 2-component \emph{links}, this is not the case (Theorem~\ref{thmMain}).   
  
\begin{theorem}\label{thmMain}
There exists a link diagram invariant  unchanged by Types I and III Reidemeister moves and can be changed by a Type II Reidemeister move.  
\end{theorem}
Here the definition of a \emph{link diagram invariant} is given by extending the notion of a \emph{knot diagram invariant} as above  to the multi-component case.  In  this paper, we prove  Theorem~\ref{thmMain} by giving an explicit Gauss diagram formula. 
\section{Proof of Theorem~\ref{thmMain}}
In this section, the definitions  \cite[Page~299, Section~2.3]{Ostlund2004} of \emph{Gauss diagram} and \emph{Gauss diagram formula},  notations (symbols)  \cite[Figures~2--4]{Polyak2010} of oriented Reidemeister moves $\Omega_{*}$ obey \cite{Ostlund2004, Polyak2010} except that  ``Gauss diagram formulas" are called ``arrow diagram formulas" in \cite{Ostlund2004}.     

\subsection{Invariance under Type I and III Reidemeister moves in a single component}
In this case, any Type I Reidemeister move belongs to a single component.   
Then it is obvious that the Gauss diagram formula $\langle \fonep, \cdot \rangle$ is unchanged by $\Omega_{1*}$ ($*=a, b, c, d$) since $\langle \fonep, \cdot \rangle$ has no isolated arrows (note also that there exists an endpoint of another arrow between the two endpoints of the arrow having the plus sign).  
It is also clear that if $\Omega_{3*}$ ($*=a, b, c, d, e, f, g, h$) is applied completely in a single component, $\langle \fonep, \cdot \rangle$ is unchanged by $\Omega_{3*}$ ($*=a, b, c, d, e, f, g, h$) since $\langle \fonep, \cdot \rangle$ has only one  arrow in a single circle.    
\subsection{The difference of values before and after applying Type II Reidemeister moves in a single   component}\label{sec:1II}
If the arrow in a single circle of $\langle \fonep, \cdot \rangle$ responds to the arrow from two crossings created by applying a single $\Omega_{2*}$ ($*=a, b, c, d$), exactly $+1$ is added to the value.      
\subsection{Invariance under Type II Reidemeister moves among two  components}\label{sec:2II}
Gauss diagram presentations of all  cases for $\Omega_{2*}$ ($*=a, b, c, d$) are given by \cite[$\Omega_{II+\pm}, \Omega_{II -\pm}$, Table~1]  {Ostlund2004}.  For $\Omega_{II+-}$, this invariance had already been  given by \cite[Page~11]{ItoOyamaguchi2020}.  By just replacing the paired signs ``$+-$" of $\Omega_{II+-}$ with each of the other cases, we have the proof of the corresponding invariance.   
\subsection{Invariance under $\Omega_{3a}$ among two  components}\label{sec:3a}
The Gauss diagram presentation of  $\Omega_{3a}$ among two  components is given by \cite[$\Omega_{III+-+*}$ ($*=b, m, t$),  Table~1]{Ostlund2004}; these invariances of $\Omega_{III+-+*}$ ($*=b, m, t$) had already been given by \cite[Page~11]{ItoOyamaguchi2020}.
\subsection{Invariance under $\Omega_{3*}$ (:any type) among two components}
First, we note that Sections~\ref{sec:1II}, \ref{sec:2II}, and \ref{sec:3a} imply  Table~\ref{table:Difference}.   
\begin{table}[h!]
\caption{``Positive" means the direction of the move increasing  crossings.}\label{table:Difference}
\begin{tabular}{|c|c|c|c|}\hline
Move  & positive $\Omega_{2*}$ in one  component & $\Omega_{2*}$ among two components & $\Omega_{3a}$ \\ \hline
Difference & $+1$ & $0$ & $0$ \\  \hline
\end{tabular}
\end{table}
$\Omega_{2*}$ $*=a, b, c, d$ is positive (negative,~resp.) if the Reidemeister move increasing (decreasing,~resp.) crossings.    
\begin{table}
\caption{Each Type III move $M$ is decomposed by the sequence consisting of three Reidemeister moves.  If Type II move among two components, the other Type II move in the same line should be the move  among two components.}\label{table:3rd}
\begin{tabular}{|c|c|c|c|c|}\hline
Move $M$ & Type II move & Type III  & Type II move & Reference  \\ \hline
$\Omega_{3b}$  & positive $\Omega_{2c}$ & $\Omega_{3a}$ & negative $\Omega_{2d}$ & \cite[Lemma~2.3]{Polyak2010}  \\ \hline
$\Omega_{3c}$  & positive $\Omega_{2c}$ & $\Omega_{3a}$ & negative $\Omega_{2d}$ & \cite[Lemma~2.4]{Polyak2010}  \\ \hline
$\Omega_{3d}$  & positive $\Omega_{2a}$ & $\Omega_{3b}$ & negative $\Omega_{2b}$ & \cite[Lemma~2.6 (1st line)]{Polyak2010}  \\ \hline
$\Omega_{3e}$  & positive $\Omega_{2a}$ & $\Omega_{3b}$ & negative $\Omega_{2b}$ & \cite[Lemma~2.6 (2nd line)]{Polyak2010}  \\ \hline
$\Omega_{3f}$  & positive $\Omega_{2d}$ & $\Omega_{3a}$ & negative $\Omega_{2c}$ & \cite[Lemma~2.6 (3rd line)]{Polyak2010}  \\ \hline
$\Omega_{3g}$  & positive $\Omega_{2c}$ & $\Omega_{3f}$ & negative $\Omega_{2d}$ & \cite[Lemma~2.6 (4th line)]{Polyak2010}  \\ \hline
$\Omega_{3h}$  & positive $\Omega_{2a}$ & $\Omega_{3g}$ & negative $\Omega_{2b}$ & \cite[Lemma~2.6 (5th line)]{Polyak2010}  \\ \hline
\end{tabular}
\end{table}
Second, we note that in Table~\ref{table:3rd} derived from \cite{Polyak2010}, if a Reidemeister move is Type II in a single component (among two components,~resp.), the other Type II move in the same line should be the move in a single component (among two components,~resp.).   This fact together with Tables~\ref{table:Difference} and \ref{table:3rd} implies the invariance under any Type III move. 
Here note that any type $\Omega_{3*}$ is given by a single $\Omega_{3a}$ and even number $2k$ of $\Omega_{2*}$ where exactly $k$ times $\Omega_{2*}$ are of type increasing crossings and the others are of type decreasing crossings. 
\subsection{Detecting the necessity of Type II Reidemeister move}
\begin{figure}[h!]
\includegraphics[width=5cm]{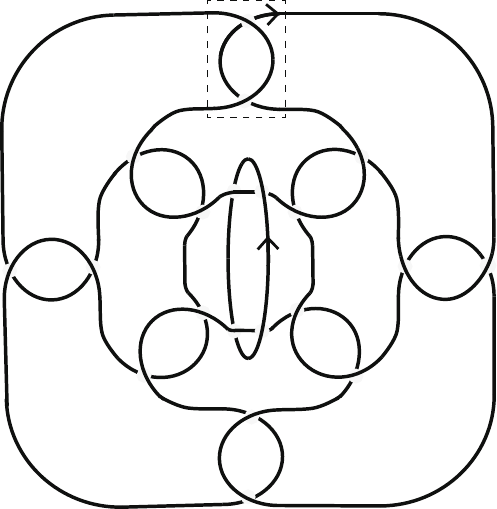}
\\
\vspace{1cm}
\includegraphics[width=8cm]{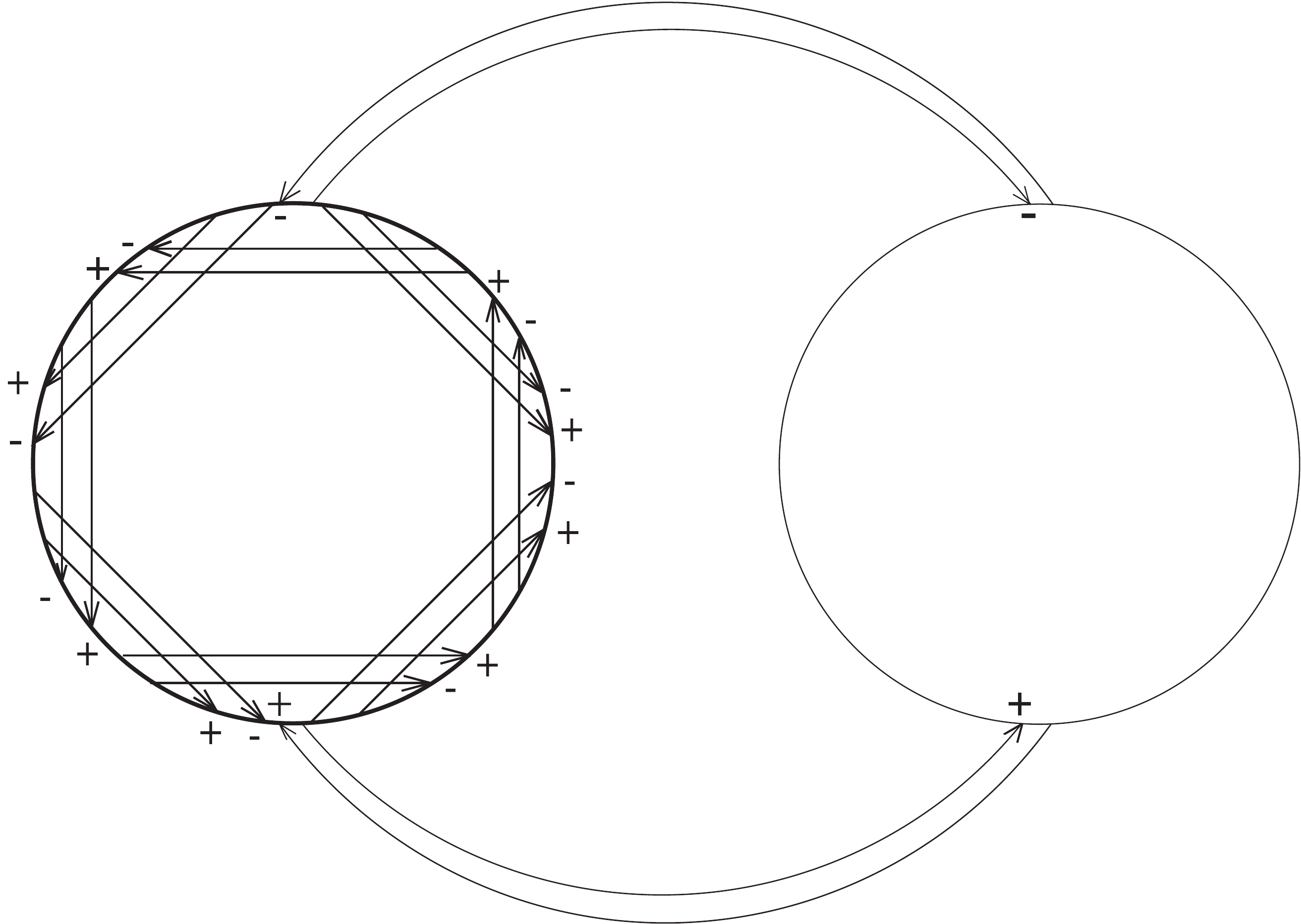}
\caption{A link diagram (upper) and its Gauss diagram (lower)}\label{fig:Tri}
\end{figure}

\begin{figure}[h!]
\includegraphics[width=8cm]{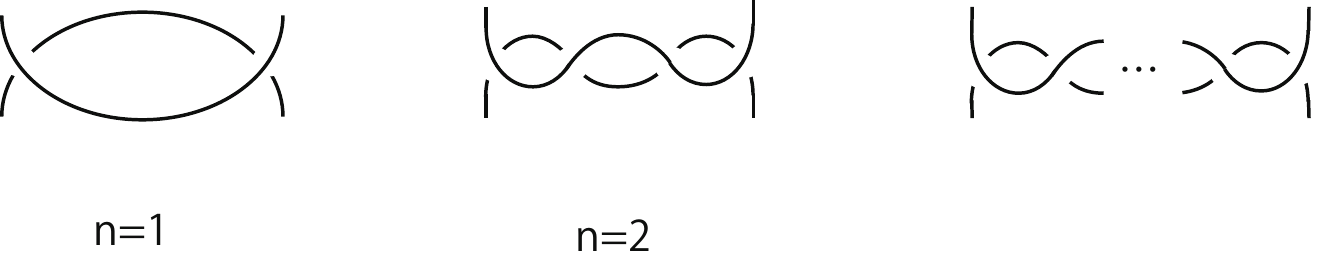}
\caption{From a single bigon given by a Type II Reidemeister moves (leftmost) to $2n-1$ bigons given by applying Type II Reidemeister moves $n$ times (rightmost); $n=2$ case (center)}\label{fig:bigon}
\end{figure}

A link diagram invariant $\langle \fonep, D_L \rangle$ detects the necessity of Type II Reidemeister move as follows (Lemma~\ref{lem:Unlink}): 
\begin{lemma}\label{lem:Unlink}
Let $D_L$ be a diagram of the trivial link as in Fig.~\ref{fig:Tri} and let $D_U$ be the 2-component diagram has no crossings.   
$\langle \fonep, D_L \rangle$ $=$ $-1$ and $\langle \fonep, D_U \rangle$ $=$ $0$.  
\end{lemma}
It is easy to generalize Lemma~\ref{lem:Unlink} to an infinitely many link diagrams.  
We focus on the single bigon (marked by the dotted box as in Fig.~\ref{fig:Tri}); we note that this bigon is obtained from a Type II Reidemeister move.  If we apply Type II Reidemeister moves $n$ times here, then we have successive $2n-1$ bigons as in Fig.~\ref{fig:bigon} (rightmost).      
Let $D_{L(n)}$ be the link diagram by replacing the single bigon with $2n-1$ bigons as in Fig.~\ref{fig:bigon} (rightmost).        
Then $\langle \fonep, D_{L(n)} \rangle$ $=$ $-n$.  


\section{Remarks}
\begin{itemize}
\item If we use $\fonem$ and the corresponding link diagram, we have  essentially the same result as the above.  
\item The above example is constructed seeing \cite{HaggeYazinski2014}.  Further, it is easy to generalize the above example to give many other examples of link diagrams using \cite{ItoTakimura2016}.  
\item The above discussion holds for virtual links, multi-component plane curves, spherical curves, and knot projections.  
\item Some important problems on this topic were treated and formulated in noughties \cite[Appendix]{Manturov2004} and \cite{Hagge2006, HaggeYazinski2014}.      Although this topic in the last decades had been treated by some works \cite{ChengGao2012, OshiroShimizuYaguchi2017} separately, any universal technique is a few \cite{ItoTakimura2016, ItoTakimura2020RII} and is encouraged.      
For example, since every knot projection is interpreted as an ascending or a descending diagram, infinitely many examples showing ``the necessity of the second Reidemeister moves of knot diagrams" had been already given in \cite{ItoTakimura2016} for virtual knots, plane curves, and spherical curves.   
\item It is observed that $n$ (appearing in the proof as above) looks the RI\!I number $n$ \cite{ItoTakimura2020RII} of knot projections.     
\end{itemize}

\section*{Acknowledgements}
The work is partially supported by JSPS KAKENHI Grant Number 20K03604.  
The author would like to thank Ms. Nanako Akiyama for giving me an electronic figure.

\bibliographystyle{plain}
\bibliography{Ref}
\end{document}